\documentclass[10pt]{article}

\usepackage{amsmath}
\usepackage{amssymb}
\usepackage{indentfirst}
\usepackage{graphics} 
\usepackage{color}

\makeatletter

\@addtoreset{equation}{section}
\makeatother
\def\<{\langle}
\def\>{\rangle}

\newtheorem{lem}{Lemma}[section]
\newtheorem{theo}{Theorem}[section]
\newtheorem{rem}{Remark}[section]
\newtheorem{pro}{Proposition}[section]

\makeatletter
   
   \@addtoreset{equation}{section}
\makeatother

\begin{document}
\title{\bf Energy decay for wave equations with\\ a potential and a localized damping}
\author{Ryo Ikehata\thanks{ikehatar@hiroshima-u.ac.jp} \\ {\small Department of Mathematics, Division of Educational Sciences}\\ {\small Graduate School of Humanities and Social Sciences} \\ {\small Hiroshima University} \\ {\small Higashi-Hiroshima 739-8524, Japan} \\and\\ Xiaoyan Li \thanks{Corresponding author: xiaoyanli@hust.edu.cn} \\ {\small School of Mathematics and Statistics}\\ {\small Huazhong University of Science and Technology} \\ {\small Wuhan, Hubei 430074, PR China}}
\date{}
\maketitle
\begin{abstract}
We consider the total energy decay together with $L^{2}$-bound of the solution itself of the Cauchy problem for wave equations with a localized damping and a short-range potential. We treat it in the one dimensional Euclidean space ${\bf R}$. We adopt a simple multiplier method to study them. In this case, it is essential that the compactness of support of the initial data is not assumed. Since this problem is treated in the whole space, the Poincar\'e and  Hardy inequalities are not available as is developed in the exterior domain case for $n \geq 1$. For compensating  such a lack of useful tools, the potential plays an effective role. As an application,  the global existence of small data solution for a semilinear problem is provided. 
\end{abstract}
\section{Introduction}
\footnote[0]{Keywords and Phrases: Wave equation; one dimensional space; potential; localized damping; multiplier method; energy decay.}
\footnote[0]{2000 Mathematics Subject Classification. Primary 35L70; Secondary 35L05, 35B33, 35B40.}

We consider the Cauchy problem
 for wave equations with a potential and a localized damping in one dimensional Euclidean space ${\bf R}$
\begin{equation}
u_{tt}(t,x) - u_{xx}(t,x) + V(x)u(t,x) + a(x)u_{t}(t,x) = 0,\ \ \ (t,x)\in (0,\infty)\times {\bf R},\label{eqn}
\end{equation}
\begin{equation}
u(0,x)= u_{0}(x),\ \ u_{t}(0,x)= u_{1}(x),\ \ \ x\in {\bf R},\label{initial}
\end{equation}
where the initial data $[u_{0},u_{1}]$ are taken from the usual energy space for the moment
\[u_{0} \in H^{1}({\bf R}),\quad u_{1} \in L^{2}({\bf R}),\]
and we denote
\[u_{t}=\frac{\partial u}{\partial t},\quad u_{tt}=\frac{\partial^2 u}{\partial t^2},\quad u_{xx}=\frac{\partial^2 u}{\partial x^2}.\]

Throughout this paper, $\| \cdot\|_q$ stands for the usual $L^q({\bf R})$-norm. For simplicity of notation, in particular, we use $\| \cdot\|$ instead of $\| \cdot\|_2$. $H^{1}({\bf R})$-norm is denoted by $\Vert\cdot\Vert_{H^{1}}$. The total energy $E_{u}(t)$ to the solution $u(t,x)$ of (1.1) is defined by
\begin{equation}
E_{u}(t)=\frac{1}{2}(\| u_t(t,\cdot)\|^2+\| u_x(t,\cdot)\|^2+\|\sqrt{V(\cdot)}u(t,\cdot)\|^2 ).
\end{equation}
$f \in {\rm BC}({\bf R})$ implies that
$f$ is continuous and bounded in ${\bf R}$, and $f \in {\rm BC}^{1}({\bf R})$ means
$f, f' \in  {\rm BC}({\bf R})$. \\

We shall impose the following two assumptions on $a(x)$:\\
({\bf A.1})\,$a \in {\rm BC}({\bf R})$ and $a(x) \geq 0$,\\
({\bf A.2})\,there exists a constant $L > 0$ and $\varepsilon_{1} > 0$ such that
\[a(x) \geq \varepsilon_{1},~~~ \vert x\vert \geq L.\]
\noindent
Additionally, one assumes the following hypothesis on $V(x)$:\\
({\bf V.1})\,$V \in {\rm BC}^{1}({\bf R})$, $V(x) > 0$ for all $x \in {\bf R}$,\\
({\bf V.2})\,$V'(x)x \leq 0$ for all $x \in {\bf R}$.\\
\begin{rem}{\rm By ({\bf V.2}), the potential $V(x)$ is monotone increase in ${\bf R^-}$ and  monotone decrease in ${\bf R^+}$.} 
\end{rem} 
Under these conditions, it is known that for each $[u_{0},u_{1}] \in H^{1}({\bf R})\times L^{2}({\bf R})$, the problem (1.1)-(1.2) admits a unique weak solution $u \in {\rm C}([0,\infty);H^{1}({\bf R})) \cap {\rm C}^{1}([0,\infty);L^{2}({\bf R}))$ (cf. \cite{ikawa}).\\

Let us mention research background on the equation (1.1). \\
In the exterior domain case, in \cite{Nakao-1} the author derives the total energy decay estimate $E_{u}(t) = {\cal O}(t^{-1})$ as $t \to \infty$ for the mixed problem (without potential terms)
\begin{equation}\label{1-3}
u_{tt}(t,x) - \Delta u(t,x) + a(x)u_{t}(t,x) = 0,\ \ \ (t,x)\in (0,\infty)\times \Omega,
\end{equation}
\begin{equation} \label{1-4}
u(0,x)= u_{0}(x),\ \ u_{t}(0,x)= u_{1}(x),\ \ \ x\in \Omega,
\end{equation}
\begin{equation} \label{1-5}
u(t,x)= 0,\ \ \ x\in \partial\Omega,\quad t > 0, 
\end{equation}
where $\Omega = {\bf R}^{n}\setminus \bar{{\cal O}} \subset {\bf R}^{n}$ is a smooth exterior domain. In \cite{Nakao-1}, the damping $a(x)$ is effective near infinity (like our assumption ({\bf A.2})) and near a part of trapping boundary of $\partial\Omega$. Soon after \cite{Nakao-1}, under additional condition on the weighted initial data, the author in \cite{I} obtained faster decay estimates 
such as $E_{u}(t) = {\cal O}(t^{-2})$ and $\Vert u(t,\cdot)\Vert = {\cal O}(t^{-1/2})$ ($t \to \infty$) to the equation \eqref{1-3}. However, in \cite{I} the obstacle ${\cal O}$ must be star-shaped relative to some point to erase the influence of  trapped rays. Since these results rely on the Poincar\'e and/or Hardy inequalities, only exterior domain case and higher dimensional case ($n \geq 2$) were treated. So, if one considers the Cauchy problems of  \eqref{1-3} in the Euclidean space ${\bf R}^{n}$, similar results to \cite{I} can be obtained only in the higher dimensional case $n \geq 3$. In this sense, to get faster energy decay like $E_{u}(t) = {\cal O}(t^{-2})$ seems completely open for the low dimensional case $n = 1,2$ to the equation \eqref{1-3}. A generalization of \cite{I} without assuming star-shaped obstacle ${\cal O}$ was deeply studied in \cite{A}. By assuming that the boundary of obstacle ${\cal O}$ admits no trapped rays of geometric optics, the generalization of \cite{Nakao-1} and \cite{I} was also discussed more in detail in \cite{D}. Decay and non-decay properties of the total energy for  \eqref{1-3} were studied with a logarithmic type (time-space dependent) damping $a(t,x)$ instead of $a(x)$  in \cite{Mo}. On the topic of energy decay for wave equation with asymptotically periodic damping, we refer the readers to \cite{J}. Another generalization of \cite{I} for \eqref{1-3} was considered on the noncompact Riemannian manifold in \cite{Z}. Additionally, in \cite{Zua} and \cite{Nakao-2}, the authors have considered the Cauchy and exterior mixed problems for the Klein -Gordon type wave equations with localized dissipations. In this case, to capture the behavior of  $L^{2}$-norm of the solution as $t \to \infty$ seems much more easy, because the corresponding energy functional itself contains the $L^2$-norm of solution. It should be mentioned that an interesting problem awareness on the failure of the Hardy inequality in the one dimensional case is studied similarly in the paper \cite{SW}. We use a help of potential to avoid such a failure of the Hardy inequality. Diffusion phenomenon in abstract form and its applications to wave equations with variable damping coefficients can be found in \cite{RTY}, there non-degenerate and bounded damping coefficients $a(x)$ are treated, however, any potential terms are not considered.\\

The purpose of this paper is to consider whether faster energy decay estimates to problem (1.1)-(1.2) can be observed or not with the help of potential $V(x)$ in the one dimensional case. The problem itself is never trivial in the sense that one has no any Hardy's and Poincar\'e's inequalities. Thus, to derive useful estimates of the solution concerning two quantities
\[\Vert u(t,\cdot)\Vert, \quad \int_{0}^{t}\int_{{\bf R}}a(x)\vert u(s,x)\vert^{2}dxds\]
is both essential parts of analysis. As one more difficult point, one must absorb the localized $L^{2}$-norm of the solution itself into the total energy in the course of proof , that is, we have to derive the following relation such that
\begin{equation}\label{ike-01}
\int_{\vert x\vert \leq L}\vert u(t,x)\vert^{2}dx \leq CE_{u}(t)
\end{equation} 
with some $C > 0$. In \cite{I}, \eqref{ike-01} can be derived with the help of Poincar\'e inequality. In our case, this can not be available  anymore. So, one borrows a role of potential to get such bounds. This causes some restrictions to the shape of $V(x)$. In a sense, the potential $V(x)$ compensates a lack of  the Poincar\'e and Hardy inequalities. This idea has its origin in the first author's recent paper \cite{I-2}. Totally, the whole space case seems much more difficult than the exterior domain case, since useful tools are less prepared. It is important how we treat the localized $L^{2}$-norm in terms of potential because the damping is not effective near origin. 

Furthermore, as is pointed out recently in \cite{Ike-JHDE}, in the case of $V(x) = a(x) = 0$ for all $x \in {\bf R}$ (i.e., free wave case), the $L^{2}$-norm of the solution itself to problem (1.1)-(1.2) grows to infinity as $t \to \infty$
\begin{equation}\label{ike-02}
\Vert u(t,\cdot)\Vert \sim \sqrt{t}\quad (t \to \infty).
\end{equation}
Thus, as is soon imaged, if non-trivial $V(x)$ and $a(x)$ are considered with a quite large $L$ (less strongly effective damping) and rapidly decaying potential $V(x)$ (less strongly effective potential), such a singular property \eqref{ike-02} naturally affects on the boundedness of $L^{2}$-norm of the solution itself. This causes some difficulty in deriving the priori-estimates of solutions. So, it seems that one dimensional case is interesting to study under small effects of damping and potential.
%

To state our results, we introduce the weighted function spaces. 

Set
$$w(x) := 1+V(x)^{-1}.$$
Then, the weighted $L^{2}$-space is defined by
$$L^{2}({\bf R}, w) := \{u \in L^{2}({\bf R})\,:\,\int_{{\bf R}}\vert u(x)\vert^{2}w(x)dx < +\infty\}$$
equipped with its norm
\[\Vert u\Vert_{L^{2}({\bf R}, w)} := \left(\int_{{\bf R}}\vert u(x)\vert^{2}w(x)dx\right)^{1/2}.\]
Note that $w^{-1}, w \in L^{1}_{loc}({\bf R})$. 

Our new result reads as follows. 
\begin{theo} \label{th1}
Suppose {\rm ({\bf A.1})}, {\rm ({\bf A.2})}, {\rm ({\bf V.1})} and {\rm ({\bf V.2})} with a fixed constant $L > 0$. Then, there exists a constant $C^{*} > 0$ such that if $V(0) < (4C^{*})^{-1}$, the weak solution $u(t,x)$ to problem {\rm (1.1)-(1.2)} with initial data $[u_{0},u_{1}] \in (H^{1}({\bf R})\cap L^{2}({\bf R},w))\times L^{2}({\bf R},w)$ satisfies
\[\Vert u(t,\cdot)\Vert \leq CI_{0},\quad E_{u}(t) \leq CI_{0}^{2}(1+t)^{-1}\]
with some constant $C > 0$, where 
\[I_0^\mu:=\|u_0\|_{H^1}^\mu + \|u_1\|^\mu +\| \frac{ u_{1}+a(\cdot)u_{0}}{\sqrt{V(\cdot)}}\|^\mu, \quad (\mu = 1~ \text{or} ~2).\]
\end{theo}
\begin{rem}
{\rm $C^{*} > 0$ is a constant closely related with the modified Poincar\'e inequality (see Lemma 2.2). For the condition $V(0) < 1/(4C^{*})$, see \eqref{9} of the text below.}
\end{rem}

\begin{rem}
{\rm The obtained decay rates are slower than that studied in \cite{M} for the one dimensional case. A singularity appeared in \eqref{ike-02} for free waves may affect on a property of the quantity $\Vert u(t,\cdot)\Vert$ (as $t \to \infty$) of the solution $u(t,x)$ to problem (1.1)-(1.2).}
\end{rem}
\noindent
{\bf Example 1.}\, Suppose $\beta > 1$. One can provide a function $V(x)$ satisfying {\rm ({\bf V.1})} and {\rm ({\bf V.2})}. Indeed, one can present such an example $V \in {\rm C}^{1}({\bf R})$ satisfying 
\[V(x) = \left\{
  \begin{array}{ll}
   \displaystyle  {\frac{2V_{0}}{L^{\beta}}- {\frac{V_{0}}{L^{2\beta}}\vert x\vert^{\beta}}},& \qquad \vert x\vert \leq L,\\[0.2cm]
   \\
   \displaystyle {V_{0}\vert x\vert^{-\beta}},&
   \qquad \vert x\vert \geq L,
   \end{array} \right. \]
where $V_0$ is a positive number. Since $V(0) = \frac{2V_{0}}{L^{\beta}}$, the smallness of $V(0)$ assumed in Theorem 1.1 can be realized by choosing small $V_{0}$ for each fixed $L > 0$. $V(x)$ is the short-range potential, and in particular, $\beta = 2$ corresponds to the scale invariant case.\\
\noindent
{\bf Example 2.}\, One can give another example by $V(x) = V_{0}e^{-\nu x^{2}}$ with $\nu > 0$ and small $V_{0}$ determined in Theorem 1.1.\\
\begin{rem}{\rm The author in \cite{NO} treats the case $a(x) = $ constant $> 0$ in ${\bf R}$, and the potential satisfies $V(x) \geq k_{0}(1+\vert x\vert)^{-\beta}$ with $0 \leq \beta < 1$ ($k_{0} > 0$), that is, the long-range potential case is considered under compact support condition on the initial data. Then, the exponential decay of the total energy is obtained. A strong role of potential is effective. While, we are treating the short-range potential together with a localized damping,  the effect of them is less strong.}
\end{rem}
\begin{rem}{\rm Our next project is to study the Cauchy problem of the equation from a similar point of view:
\[u_{tt}-\Delta u + V(x)u + a(x) u_{t} = 0,\]
where $V(x) = {\cal O}(\vert x\vert^{-\alpha})$ and $a(x) = {\cal O}(\vert x\vert^{-\beta})$ as $\vert x\vert \to \infty$ for some $\alpha > 0$, $\beta > 0$. We treat the equation in the one dimensional whole space. This seems still open (cf. \cite{G}, \cite{lai}).}
\end{rem}

This paper is organized as follows. In Section 2, we shall prove Theorem 1.1 by relying on a sophisticated multiplier method which was introduced by the first author. An application to semilinear problem of (1.1) will be presented in Section 3.

\section{Proof of Theorem 1.1}

In this section, we prove Theorem 1.1 by dividing the proof into several lemmas.
 
Firstly, one prepares the following important lemma, which plays an alternative role of the Poincar\'e inequality.
\begin{lem}\label{lemma 2.1}\,Set $V_{L} := \min\{V(L),V(-L)\} > 0$. Suppose {\bf (V.1)} and {\bf (V.2)}. Then, it holds that
\[\int_{\vert x\vert \leq L}\vert u(t,x)\vert^{2}dx \leq \frac{2}{V_{L}}E_{u}(t), \quad (t \geq 0),\]
where $u(t,x)$ is the solution to problem {\rm (1.1)-(1.2)}.
\end{lem}
\underline{{\it Proof of Lemma 2.1.}} Because of {\bf (V.1)} and {\bf (V.2)}, once one notices the relation $V_{L} \leq V(x)$ for $\vert x\vert \leq L$, it is easy to derive the following inequality
\begin{align}
\int_{\vert x\vert \leq L}\vert u(t,x)\vert^{2}dx &= \int_{\vert x\vert \leq L} 
\frac{1}{V_{L}}V_{L}\vert u(t,x)\vert^{2}dx \notag \\
&\leq \frac{2}{V_{L}}\int_{{\bf R}}\frac{1}{2}V(x)\vert u(t,x)\vert^{2}dx\notag \\
& \leq \frac{2}{V_{L}}E_{u}(t),
\end{align}
where the definition of $E_{u}(t)$ is used. 
\hfill
$\Box$
\vspace{0.2cm}

We additionally prepare the Poincar\'e type inequality in the one dimensional whole space. An essential part of its proof is found in \cite[Lemma 2.1]{mochi-2}.
\begin{lem}\label{lemma 2.2} Let $L > 0$ be a constant. Then, there is a constant $C^{*} > 0$ which depends  on $L$ such that
\[\int_{\vert x\vert \leq L}\vert w(x)\vert^{2}dx \leq C^{*}\left(\int_{{\bf R}}\vert w_{x}(x)\vert^{2}dx + \int_{\vert x\vert \geq L}\vert w(x)\vert^{2}dx\right)\]
for $w \in H^{1}({\bf R})$.
\end{lem}
\vspace{0.5cm}
The first part of the proof of Theorem 1.1 follows \cite{Nakao-1}. For this trial, one firstly considers the smooth initial data case   $[u_{0},u_{1}] \in C_{0}^{\infty}({\bf R})\times C_{0}^{\infty}({\bf R})$. Then, the corresponding solution $u(t,x)$ to problem (1.1)-(1.2) becomes sufficiently smooth to guarantee the integration by parts.\\

We also prepare the following identities.
\begin{pro}\,Let $u(t,x)$ be a smooth solution to problem {\rm (1.1)}-{\rm (1.2)} with smooth initial data $[u_{0},u_{1}] \in C_{0}^{\infty}({\bf R})\times C_{0}^{\infty}({\bf R})$. Then, it holds that
\begin{equation}\label{id1}
E_{u}(t) + \int_{0}^{t}\int_{{\bf R}}a(x)\vert u_{s}(s,x)\vert^{2}dxds = E_{u}(0),
\end{equation}
\begin{equation}\label{id2}
\frac{d}{dt}(u_{t}(t,\cdot),u(t,\cdot)) - \Vert u_{t}(t,\cdot)\Vert^{2} + \Vert u_{x}(t,\cdot)\Vert^{2} + \Vert\sqrt{V(\cdot)}u(t,\cdot)\Vert^{2} + \frac{1}{2}\frac{d}{dt}\int_{{\bf R}}a(x)\vert u(t,x)\vert^{2}dx = 0.
\end{equation}
\end{pro}

Take
\[\phi(x) = \left\{
  \begin{array}{ll}
   \displaystyle{\varepsilon_{1}},&
       \qquad \vert x\vert \leq L, \\[0.2cm]
   \displaystyle{\frac{L\varepsilon_{1}}{r}},& \qquad \vert x\vert \geq L,
   \end{array} \right. \]
where $r=|x|$. Note that $\phi(x)$ is Lipschitz continuous in ${\bf R}$.   
As in \cite{Nakao-1}, multiplying both sides of (1.1) by $\phi(x)xu_{x}$ and integrating over ${\bf R}$, by the integration by parts, one finds that
\begin{align}
&\frac{d}{dt}\int_{{\bf R}}u_{t}(t,x)\phi(x)(x\cdot u_{x}(t,x))dx + \frac{1}{2}\int_{{\bf R}}(\phi(x)+\phi'(x)x)\vert u_{t}(t,x)\vert^{2}dx \notag\\
&+ \frac{1}{2}\int_{{\bf R}}(\phi(x)+\phi'(x)x)\vert u_{x}(t,x)\vert^{2}dx
- \frac{1}{2}\int_{{\bf R}}(V(x)\phi(x)+V(x)\phi'(x)x)\vert u(t,x)\vert^{2}dx\notag\\
 &- \frac{1}{2}\int_{{\bf R}}V'(x)\phi(x)x\vert u(t,x)\vert^{2}dx
+\int_{{\bf R}}a(x)u_{t}(t,x)\phi(x)x\cdot u_{x}(t,x)dx = 0,\label{3}
\end{align}
where $\phi'(x) = \frac{d\phi}{dx}$, and $V'(x) = \frac{dV}{dx}$. Thus, it follows from \eqref{id1}, \eqref{id2} and \eqref{3} that
\begin{align}\label{4}
\frac{d}{dt}&\left(\int_{{\bf R}}u_{t}(t,x)\phi(x)(x\cdot u_{x}(t,x))dx + \alpha(u_{t}(t,\cdot),u(t,\cdot)) + \frac{\alpha}{2}\int_{{\bf R}}a(x)\vert u(t,x)\vert^{2}dx + kE_{u}(t)\right) \notag\\
+&\int_{{\bf R}}(\frac{\phi(x)+x\phi'(x)}{2}-\alpha + \frac{ka(x)}{2}+\frac{ka(x)}{2})\vert u_{t}(t,x)\vert^{2}dx \notag \\ 
+& \int_{{\bf R}}(\alpha +\frac{\phi(x)+\phi'(x)x}{2})\vert u_{x}(t,x)\vert^{2}dx + \frac{1}{2}\int_{{\bf R}}V(x)\left(2\alpha-\phi'(x)x\right)\vert u(t,x)\vert^{2}dx \notag \\ 
= -&\int_{{\bf R}}a(x)u_{t}(t,x)\phi(x)(x\cdot u_{x}(t,x))dx + \frac{1}{2}\int_{{\bf R}}\left(V'(x)x+V(x)\right)\phi(x)\vert u(t,x)\vert^{2}dx,
\end{align}
where one has just used two identities \eqref{id1} and \eqref{id2} multiplied by positive parameters $k > 0$ and $\alpha > 0$, respectively. 

Now, by {\bf (V.2)} and Lemma \ref{lemma 2.2}, it follows that
\begin{align}
	&~~~\int_{{\bf R}}\left(V'(x)x+V(x)\right)\phi(x)\vert u(t,x)\vert^{2}dx \notag\\
		&\leq \int_{\bf R}V(x)\phi(x)\vert u(t,x)\vert^{2}dx \notag \\
	&\leq \int_{\vert x\vert \leq L}V(x) \varepsilon_{1} \vert u(t,x)\vert^{2}dx+\int_{\vert x\vert \geq L} \frac{L \varepsilon_1}{|x|}V(x)\vert u(t,x)\vert^{2}dx  \notag \\
& \leq \varepsilon_{1}V(0)C^{*}\left( \int_{\vert x\vert \geq L}\vert u(t,x)\vert^{2}dx +  \int_{{\bf R}}\vert u_{x}(t,x)\vert^{2}dx \right)
+\varepsilon_1 V_L^{'} \int_{{|x|\geq L}} |u(t,x)|^2 dx \notag\\
&= \varepsilon_{1}V(0)C^{*}\int_{{\bf R}}\vert u_{x}(t,x)\vert^{2}dx + \big(V(0)C^{*}+V_L{'}\big)\int_{\vert x\vert \geq L}\varepsilon_{1}\vert u(t,x)\vert^{2}dx,
\end{align}
where $V_L{'}:= \max\{V(L),V(-L)\}$. By using {\bf (A.2)}, one has
\begin{align}
&\int_{{\bf R}}\left(V'(x)x+V(x)\right)\phi(x)\vert u(t,x)\vert^{2}dx \notag\\
\leq & \varepsilon_{1}V(0)C^{*}\int_{{\bf R}}\vert u_{x}(t,x)\vert^{2}dx + \big(V(0)C^{*}+V_L{'}\big)\int_{{\bf R}}a(x)\vert u(t,x)\vert^{2}dx.\label{5}
\end{align}

Therefore, \eqref{4} and \eqref{5} yield
\begin{align}\label{6}
\frac{d}{dt}&\left(\int_{{\bf R}}u_{t}(t,x)\phi(x)(x\cdot u_{x}(t,x))dx + \alpha(u_{t}(t,\cdot),u(t,\cdot)) + \frac{\alpha}{2}\int_{{\bf R}}a(x)\vert u(t,x)\vert^{2}dx + kE_{u}(t)\right)\notag\\
&+\int_{{\bf R}}\left(\frac{\phi(x)+x\phi'(x)}{2}-\alpha + \frac{ka(x)}{2}\right)\vert u_{t}(t,x)\vert^{2}dx + \frac{k}{2}\int_{{\bf R}}a(x)\vert u_{t}(t,x)\vert^{2}dx\notag\\
&+ \int_{{\bf R}}(\alpha + \frac{\phi(x) + \phi'(x)x}{2}-\frac{C^{*}\varepsilon_{1}V(0)}{2})\vert u_{x}(t,x)\vert^{2}dx + \frac{1}{2}\int_{{\bf R}}V(x)(2\alpha-\phi'(x)x)\vert u(t,x)\vert^{2}dx\notag\\
\leq& -\int_{{\bf R}}a(x)u_{t}(t,x)\phi(x)(x\cdot u_{x}(t,x))dx + \frac{V(0)C^{*}+V_L{'}}{2}\int_{{\bf R}}a(x)\vert u(t,x)\vert^{2}dx.
\end{align}

Now, according to  {\bf (A.2)}, there exists a  number $\alpha > 0$ and a small number $\varepsilon_{2} > 0$ such that
\begin{equation}\label{7}
\frac{\phi(x)+x\phi'(x)}{2}-\alpha + \frac{ka(x)}{2} > \varepsilon_{2} > 0,
\end{equation}
\begin{equation}\label{8}
\alpha+\frac{\phi(x)+x\phi'(x)}{2} > \varepsilon_{2} > 0
\end{equation}
with some constant $\varepsilon_{2} \in (0,\displaystyle{\frac{\varepsilon_{1}}{2}})$ and $k \geq 2$. In fact, we can choose
\[\alpha = \frac{\varepsilon_{1}}{4},\quad \varepsilon_{2} = \frac{\varepsilon_{1}}{8}.\] 
Furthermore, choose $V(0) > 0$ small enough to satisfy
\[\varepsilon_{2} = \frac{\varepsilon_{1}}{8} > \frac{C^{*}\varepsilon_{1}V(0)}{2},\]
that is, 
\begin{equation}\label{9}
0 < V(0) < \frac{1}{4C^{*}}.
\end{equation}
It's obvious that \eqref{8} and \eqref{9} guarantee the positiveness such that 
\[\gamma_{0} := 2(\varepsilon_{2}-\frac{C^{*}\varepsilon_{1}V(0)}{2}) > 0.\]
\noindent
One hand, one can obtain the following estimate
\begin{equation}\label{10}
J_{1}(t) := \left\vert\int_{{\bf R}}a(x)u_{t}\phi(x)(x\cdot u_{x})dx\right\vert\ \leq \frac{L\varepsilon_{1}k}{8}\int_{{\bf R}}a(x)\vert u_{t}\vert^{2}dx + \frac{2L\varepsilon_{1}\Vert a\Vert_{\infty}}{k}\int_{{\bf R}}\vert u_{x}\vert^{2}dx.
\end{equation} 
Indeed,   
\[J_{1}(t) \leq L \int_{\vert x\vert \leq L}\frac{\sqrt{k}}{2}\sqrt{a(x)}\vert u_{t}(t,x)\vert\varepsilon_{1}\cdot\sqrt{a(x)}\vert u_{x}(t,x)\vert\frac{2}{\sqrt{k}}dx~~~~~~~~~~~~\]
\[+ \int_{\vert x \vert \geq L}\frac{\sqrt{k}}{2}\sqrt{a(x)}\vert u_{t}(t,x)\vert\varepsilon_{1}L\sqrt{a(x)}\vert u_{x}(t,x)\vert\frac{2}{\sqrt{k}}dx\]
\[= L\varepsilon_{1}\int_{{\bf R}}\left(\frac{\sqrt{k}}{2}\sqrt{a(x)}\vert u_{t}(t,x)\vert\right)\left(\sqrt{a(x)}\vert u_{x}(t,x)\vert\frac{2}{\sqrt{k}}\right)dx\]
\[\leq \frac{L\varepsilon_{1}}{2}\int_{{\bf R}}\frac{k}{4}a(x)\vert u_{t}(t,x)\vert^{2}dx + \frac{L\varepsilon_{1}}{2}\int_{{\bf R}}\frac{4}{k}a(x)\vert u_{x}(t,x)\vert^{2}dx\]
\[\leq \frac{L\varepsilon_{1}k}{8}\int_{{\bf R}}a(x)\vert u_{t}(t,x)\vert^{2}dx + \frac{2L\varepsilon_{1}}{k}\Vert a\Vert_{\infty}\int_{{\bf R}}\vert u_{x}(t,x)\vert^{2}dx,\]
which implies the desired estimate. So, setting
\begin{equation}\label{11}
G_{k}(t) := \int_{{\bf R}}u_{t}(t,x)\phi(x)(x\cdot u_{x}(t,x))dx + \alpha(u_{t}(t,\cdot),u(t,\cdot)) + \frac{\alpha}{2}\int_{{\bf R}}a(x)\vert u(t,x)\vert^{2}dx + kE_{u}(t),
\end{equation}
from \eqref{6}, \eqref{7}, \eqref{8}, \eqref{9} and \eqref{10}, one obtains
\begin{align}
&\frac{d}{dt}G_{k}(t) + \varepsilon_{2}\int_{{\bf R}}\vert u_{t}(t,x)\vert^{2}dx + 2(\varepsilon_{2}-\frac{C^{*}\varepsilon_{1}V(0)}{2})\int_{{\bf R}}\frac{1}{2}\vert u_{x}(t,x)\vert^{2}dx \notag\\
&+ \frac{1}{2}\int_{{\bf R}}V(x)(2\alpha-\phi'(x)x)\vert u(t,x)\vert^{2}dx + \frac{k}{2}\int_{{\bf R}}a(x)\vert u_{t}(t,x)\vert^{2}dx\notag\\
\leq& \frac{L\varepsilon_{1}k}{8}\int_{{\bf R}}a(x)\vert u_{t}(t,x)\vert^{2}dx + \frac{4L\varepsilon_{1}\Vert a\Vert_{\infty}}{k}\int_{{\bf R}}\frac{1}{2}\vert u_{x}(t,x)\vert^{2}dx + \frac{V(0)C^{*}+V_L{'}}{2}\int_{{\bf R}}a(x)\vert u(t,x)\vert^{2}dx \notag,
\end{align}
which implies
\begin{align}
&\frac{d}{dt}G_{k}(t) + \frac{\varepsilon_{1}}{4}\int_{{\bf R}}\frac{1}{2}\vert u_{t}(t,x)\vert^{2}dx + (\gamma_{0}- \frac{4L\varepsilon_{1}\Vert a\Vert_{\infty}}{k})\int_{{\bf R}}\frac{1}{2}\vert u_{x}(t,x)\vert^{2}dx\notag\\
&\quad+ \frac{1}{2}\int_{{\bf R}}V(x)(2\alpha-\phi'(x)x)\vert u(t,x)\vert^{2}dx\notag\\
\leq& \frac{L\varepsilon_{1}k}{8}\int_{{\bf R}}a(x)\vert u_{t}(t,x)\vert^{2}dx +   \frac{V(0)C^{*}}{2}\int_{{\bf R}}a(x)\vert u(t,x)\vert^{2}dx\notag\\
=  &-\frac{L\varepsilon_{1}k}{8}E_{u}'(t) +  \frac{V(0)C^{*}+V_L{'}}{2}\int_{{\bf R}}a(x)\vert u(t,x)\vert^{2}dx\notag. 
\end{align}
Choose $k > 0$ large enough such that
\[P_{0} := \gamma_{0}-\frac{4L\varepsilon_{1}\Vert a\Vert_{\infty}}{k} > 0.\]
Then, since $-\phi'(x)x \geq 0$ (a.e. $x \in {\bf R}$), it follows that
\[\frac{d}{dt}G_{k}(t) + \frac{\varepsilon_{1}}{4}\int_{{\bf R}}\frac{1}{2}\vert u_{t}(t,x)\vert^{2}dx + P_{0}\int_{{\bf R}}\frac{1}{2}\vert u_{x}(t,x)\vert^{2}dx + 2\alpha\int_{{\bf R}}\frac{1}{2}V(x)\vert u(t,x)\vert^{2}dx\]
\[\leq -\frac{L\varepsilon_{1}k}{8}E_{u}'(t) +  \frac{V(0)C^{*}+V_L{'}}{2}\int_{{\bf R}}a(x)\vert u(t,x)\vert^{2}dx.\]
\noindent
Therefore, we find that there exists a positive constant 
$$\eta_{0}:=\min\{  \frac{\varepsilon_{1}}{4},~P_0,~2\alpha \}$$
depending only on $\varepsilon_{1}$, $L > 0$ and $\Vert a\Vert_{\infty}$ such that
\begin{equation}\label{12}
\frac{d}{dt}G_{k}(t) + \eta_{0}E_{u}(t) \leq -\frac{L\varepsilon_{1}k}{8}E_{u}'(t) + \frac{V(0)C^{*}+V_L{'}}{2}\int_{{\bf R}}a(x)\vert u(t,x)\vert^{2}dx.
\end{equation}

\noindent
Integrating \eqref{12} over $[0,t]$, one has
\begin{align}\label{13}
&G_{k}(t) + \eta_{0}\int_{0}^{t}E_{u}(s)ds \notag \\
\leq &G_{k}(0) - \frac{L\varepsilon_{1}k}{8}E_{u}(t) + \frac{L\varepsilon_{1}k}{8}E_{u}(0) + \frac{V(0)C^{*}+V_L{'}}{2}\int_{0}^{t}\int_{{\bf R}}a(x)\vert u(s,x)\vert^{2}dxds\notag \\
\leq &G_{k}(0) + \frac{L\varepsilon_{1}k}{8}E_{u}(0) + \frac{V(0)C^{*}+V_L{'}}{2}\int_{0}^{t}\int_{{\bf R}}a(x)\vert u(s,x)\vert^{2}dxds.   
\end{align}
Here, we note that
\begin{equation}\label{14}
\vert G_{k}(0)\vert \leq C(\Vert u_{0}\Vert_{H^{1}}^{2} + \Vert u_{1}\Vert^{2})
\end{equation}
with some $C > 0$. 
Thus, one can arrive at the following lemma.   
\begin{lem}\label{lemma 2.3}\, Let $[u_{0},u_{1}] \in C_{0}^{\infty}({\bf R})\times C_{0}^{\infty}({\bf R})$. Then, for the corresponding smooth solution $u(t,x)$ to problem {\rm (1.1)-(1.2)}, it holds that
\[G_{k}(t) + \eta_{0}\int_{0}^{t}E_{u}(s)ds \leq C(\Vert u_{0}\Vert_{H^{1}}^{2} + \Vert u_{1}\Vert^{2} + \int_{0}^{t}\int_{{\bf R}}a(x)\vert u(s,x)\vert^{2}dxds)\quad (t \geq 0),\]
with some generous constants $\eta_{0} > 0$ and $C > 0$, provided that $V(0) > 0$ and $k \geq 2$ are chosen small and large, respectively.
\end{lem}

In the next part, let us check the positiveness of $G_{k}(t)$ defined by \eqref{11}. The proof is similar to \cite[Lemma 2.3]{I} except for using Lemma \ref{lemma 2.1} instead of the Poincar\'e inequality. Indeed, for $\varepsilon > 0$, it follows from Lemma \ref{lemma 2.1} and {\bf (A.2)} that
\begin{align}\label{15}
-\alpha(u_{t}(t,\cdot), u(t,\cdot)) &\leq \frac{\alpha}{2\varepsilon}\Vert u_{t}(t,\cdot)\Vert^{2} + \frac{\alpha\varepsilon}{2}\Vert u(t,\cdot)\Vert^{2} \notag \\ 
&\leq \frac{\alpha}{\varepsilon}E_{u}(t) + \frac{\alpha\varepsilon}{2}\left(\frac{1}{\varepsilon_{1}}\int_{\vert x\vert \geq L}a(x)\vert u(t,x)\vert^{2}dx + \int_{\vert x\vert \leq L}\vert u(t,x)\vert^{2}dx\right)\notag \\ 
&\leq \frac{\alpha}{\varepsilon}E_{u}(t) + \frac{\alpha\varepsilon}{2}\left(\frac{1}{\varepsilon_{1}}\int_{\vert x\vert \geq L}a(x)\vert u(t,x)\vert^{2}dx + \frac{2}{V_{L}}E_{u}(t)\right)\notag \\ 
&\leq \left(\frac{\alpha}{\varepsilon} + \frac{\alpha\varepsilon}{V_{L}}\right)E_{u}(t) + \frac{\alpha\varepsilon}{2\varepsilon_{1}}\int_{{\bf R}}a(x)\vert u(t,x)\vert^{2}dx.
\end{align}
Furthermore, from the definition of the function $\phi(x)$ one sees
\begin{align}\label{16}
&-\int_{{\bf R}}u_{t}(t,x)\phi(x)x u_{x}(t,x)dx \notag \\
\leq & \int_{\vert x\vert \geq L}\vert u_{t}(t,x)\vert\phi(x)\vert x\vert \vert u_{x}(t,x)\vert dx + \int_{\vert x\vert \leq L}\vert u_{t}(t,x)\vert\phi(x)\vert x\vert \vert u_{x}(t,x)\vert dx \notag \\ 
\leq &\int_{\vert x\vert \geq L}\vert u_{t}(t,x)\vert\frac{L\varepsilon_{1}}{\vert x\vert}\vert x\vert \vert u_{x}(t,x)\vert dx + \int_{\vert x\vert \leq L}\vert u_{t}(t,x)\vert\varepsilon_{1}L\vert u_{x}(t,x)\vert dx\notag \\ 
=& L\varepsilon_{1}\int_{\vert x\vert \geq L}\vert u_{t}(t,x)\vert\vert u_{x}(t,x)\vert dx + \varepsilon_{1}L\int_{\vert x\vert \leq L}\vert u_{t}(t,x)\vert\vert u_{x}(t,x)\vert dx \notag \\ 
\leq &L\varepsilon_{1}\int_{{\bf R}}\vert u_{t}(t,x)\vert\vert u_{x}(t,x)\vert dx\notag \\
 \leq & L\varepsilon_{1}E_{u}(t).
\end{align}
Thus, \eqref{15} and \eqref{16} imply
\begin{align*}
&-\alpha(u_{t}(t,\cdot), u(t,\cdot)) - \int_{{\bf R}}u_{t}(t,x)\phi(x)x u_{x}(t,x)dx\\
	\leq& \left(\frac{\alpha}{\varepsilon} + \frac{\alpha\varepsilon}{V_{L}} + L\varepsilon_{1}\right)E_{u}(t) + \frac{\alpha\varepsilon}{2\varepsilon_{1}}\int_{{\bf R}}a(x)\vert u(t,x)\vert^{2}dx.
\end{align*}

Finally, taking $0<\varepsilon <\varepsilon_1$  and choosing $k \geq 2$ large enough such that 
\[\frac{\alpha}{\varepsilon} + \frac{\alpha\varepsilon}{V_{L}} + L\varepsilon_{1} < k,\] 
one can arrive at the following lemma.
\begin{lem}\label{lemma 2.4}\, Let $[u_{0},u_{1}] \in C_{0}^{\infty}({\bf R})\times C_{0}^{\infty}({\bf R})$. Then, for the corresponding smooth solution $u(t,x)$ to problem {\rm (1.1)-(1.2)}, it holds that 
\[G_{k}(t) \geq 0 \quad (t \geq 0)\]
for large $k \gg 1$.
\end{lem}

Now, one prepares the crucial result to derive main estimates of Theorem 1.1. The idea is an application of the method recently developed in \cite{I-2}.
\begin{lem}\label{lemma 2.5}\, Let $[u_{0},u_{1}] \in C_{0}^{\infty}({\bf R})\times C_{0}^{\infty}({\bf R})$. Then, for the corresponding smooth solution $u(t,x)$ to problem {\rm (1.1)-(1.2)}, it holds that
\[\Vert u(t,\cdot)\Vert^{2} + \int_{0}^{t}\int_{{\bf R}}a(x)\vert u(s,x)\vert^{2}dxds \leq  C\left(\Vert u_{0}\Vert^{2} + \int_{{\bf R}}\frac{\vert u_{1}(x)+a(x)u_{0}(x)\vert^{2}}{V(x)}dx\right) \quad (t \geq 0),\]
where $C > 0$ is a generous constant.
\end{lem} 

\underline{{\it Proof of Lemma \ref{lemma 2.5}.}}\, At first, for the solution $u(t,x)$ to problem (1.1)-(1.2), one sets 
\[v(t,x) := \int_{0}^{t}u(s,x)ds.\] 
This simple idea comes from \cite{IM}, which is a modification of the celebrated Morawetz method. It can be seen that $v(t,x)$ satisfies
\begin{equation}\label{17}
v_{tt}(t,x) - v_{xx}(t,x) + V(x)v(t,x) + a(x)v_{t}(t,x) = u_{1}(x) + a(x)u_{0}(x), \quad t > 0,\quad x \in {\bf R},
\end{equation}
\begin{equation}\label{18}
v(0,x)= 0,\hspace{0,50cm}v_{t}(0,x) = u_{0}(x), \quad x \in {\bf R}.
\end{equation} 
Multiplying both sides of \eqref{17} by $v_{t}$ and integrating it over $[0,t]\times {\bf R}$, it follows from \eqref{18} that
\[\frac{1}{2}\Vert v_{t}(t,\cdot)\Vert^{2} + \frac{1}{2}\Vert v_{x}(t,\cdot)\Vert^{2} + \frac{1}{2}\int_{{\bf R}}V(x)\vert v(t,x)\vert^{2}dx + \int_{0}^{t}\int_{{\bf R}}a(x)\vert v_{s}(s,x)\vert^{2}dxds\]
\begin{equation}\label{45}
= \frac{1}{2}\Vert u_{0}\Vert^{2} + (u_{1} + a(\cdot)u_{0}, v(t,\cdot)), \quad t \geq 0.~~~~~~~~~~~~~~~~~~~~~~~~~~~~~~~~
~~~~~~~~~~~~~~
\end{equation}
Now, let us estimate the final term of the right hand side of \eqref{45} in order to absorb it into the left hand side. Indeed, by the Schwarz inequality, one has
\begin{align}\label{20}
\vert (u_{1}+a(\cdot)u_{0}, v(t,\cdot))\vert& \leq \int_{{\bf R}}\vert u_{1}(x)+a(x)u_{0}(x)\vert\vert v(t,x)\vert dx \notag\\
&= \int_{{\bf R}}\frac{\vert u_{1}(x)+a(x)u_{0}(x)\vert}{\sqrt{V(x)}}\left(\sqrt{V(x)}\vert\vert v(t,x)\vert\right)dx\notag\\
&\leq \left(\int_{{\bf R}}\frac{\vert u_{1}(x)+a(x)u_{0}(x)\vert^{2}}{V(x)}dx \right)^{1/2}\left(\int_{{\bf R}}V(x)\vert v(t,x)\vert^{2}dx\right)^{1/2}\notag\\
&\leq \int_{{\bf R}}\frac{\vert u_{1}(x)+a(x)u_{0}(x)\vert^{2}}{V(x)}dx + \frac{1}{4}\int_{{\bf R}}V(x)\vert v(t,x)\vert^{2}dx.
\end{align}
Thus \eqref{45} and \eqref{20} imply the desired estimate
\[\frac{1}{2}\Vert u(t,\cdot)\Vert^{2} + \frac{1}{2}\Vert v_{x}(t,\cdot)\Vert^{2} + \frac{1}{4}\int_{{\bf R}}V(x)\vert v(t,x)\vert^{2}dx + \int_{0}^{t}\int_{{\bf R}}a(x)\vert u(s,x)\vert^{2}dxds\]
\[\leq \frac{1}{2}\Vert u_{0}\Vert^{2} + \int_{{\bf R}}\frac{\vert u_{1}(x)+a(x)u_{0}(x)\vert^{2}}{V(x)}dx,~~~~~~~~~~~~~~~~~~~~~~~~~~~~~~~~~~~~~~~~~~~~~~~~~~~~~~\]
because of $v_{t} = u$. 
\hfill
$\Box$\\

Lemmas \ref{lemma 2.3}, \ref{lemma 2.4} and \ref{lemma 2.5} imply the following decay estimates of the total energy. The proof is standard (cf. \cite[Lemma 2.4]{I}).

\begin{pro}\label{pro2.1}\,Let $[u_{0},u_{1}] \in C_{0}^{\infty}({\bf R})\times C_{0}^{\infty}({\bf R})$. Then, for the corresponding smooth solution $u(t,x)$ to problem {\rm (1.1)-(1.2)}, it holds that
\[\int_{0}^{t}E_{u}(s)ds \leq C\left(\Vert u_{0}\Vert_{H^{1}}^{2} + \Vert u_{1}\Vert^{2} + \int_{{\bf R}}\frac{\vert u_{1}(x)+a(x)u_{0}(x)\vert^{2}}{V(x)}dx\right)\quad (t \geq 0),\]
\[(1+t)E_{u}(t) \leq C\left(\Vert u_{0}\Vert_{H^{1}}^{2} + \Vert u_{1}\Vert^{2} + \int_{{\bf R}}\frac{\vert u_{1}(x)+a(x)u_{0}(x)\vert^{2}}{V(x)}dx\right)\quad (t \geq 0),\]
where $C > 0$ is a generous constant.
\end{pro}

As a consequence of Proposition \ref{pro2.1}, one can get the local energy decay result.
\begin{pro}\label{pro2.2}\,Let $[u_{0},u_{1}] \in C_{0}^{\infty}({\bf R})\times C_{0}^{\infty}({\bf R})$. Then, for the corresponding smooth solution $u(t,x)$ to problem {\rm (1.1)-(1.2)}, it holds that
\[(1+t)\int_{\vert x\vert \leq L}\vert u(t,x)\vert^{2}dx \leq C I_{0}^{2}\quad (t \geq 0),\]
where $C > 0$ is a generous constant.
\end{pro}
{\it Proof.} Indeed, it follows from Lemma \ref{lemma 2.1} and Proposition \ref{pro2.1} that
\[(1+t)\int_{\vert x\vert \leq L}\vert u(t,x)\vert^{2}dx \leq \frac{2}{V_{L}}(1+t)E_{u}(t) \leq CI_{0}^{2}.\]
This implies the desired estimate.
\hfill
$\Box$

Finally let us prove Theorem 1.1 with above preparations. \\
\vspace{0.2cm}
{\it \underline{Proof of Theorem 1.1.}}
Theorem 1.1 is a direct consequence of Proposition \ref{pro2.1}, Lemma \ref{lemma 2.5} and density argument (i.e., cut-off technique and the mollifier method). In fact, one can choose a sequence $[\phi_{n}, \psi_{n}] \in C_{0}^{\infty}({\bf R})\times C_{0}^{\infty}({\bf R})$ such that
\[\Vert \phi_{n}-u_{0}\Vert_{L^{2}({\bf R},w)} + \Vert \phi_{n}'-u_{0}'\Vert \to 0 \quad (n \to \infty),\]
\[\Vert \psi_{n}-u_{1}\Vert_{L^{2}({\bf R},w)} \to 0 \quad (n \to \infty),\]
where
\[w(x) := 1+V(x)^{-1}.\]
Let $u^{(n)}(t,x)$ be a corresponding smooth solution to problem (1.1)-(1.2) with initial data $u_{0} := \phi_{n}$ and $u_{1} := \psi_{n}$. It is easy to obtain the following relations between $u^{(n)}(t,x)$ and $u(t,x)$
\[\sup_{t \in [0,\infty)}\left(\Vert u_{t}^{(n)}(t,\cdot) - u_{t}(t,\cdot)\Vert + \Vert u_{x}^{(n)}(t,\cdot) - u_{x}(t,\cdot)\Vert + \Vert\sqrt{V(\cdot)}(u^{(n)}(t,\cdot) - u(t,\cdot))\Vert\right) \to 0 \quad (n \to \infty),\]  
\[\sup_{t \in [0,T]}\Vert u^{(n)}(t,\cdot) - u(t,\cdot)\Vert \to 0 \quad (n \to \infty)\] 
for each $T > 0$. Then, it follows from Lemma \ref{lemma 2.5} and Proposition \ref{pro2.1} that  
\[(1+t)E_{u^{(n)}}(t) \leq C\left(\Vert \phi_{n}\Vert_{H^{1}}^{2} + \Vert \psi_{n}\Vert^{2} + \int_{{\bf R}}\frac{\vert \psi_{n}(x)+a(x)\phi_{n}(x)\vert^{2}}{V(x)}dx \right),\quad (t \geq 0),\]
\[\Vert u^{(n)}(t,\cdot)\Vert^{2} \leq C\left(\Vert \phi_{n}\Vert^{2} + \int_{{\bf R}}\frac{\vert \psi_{n}(x)+a(x)\phi_{n}(x)\vert^{2}}{V(x)}dx\right),\quad (t \geq 0).\]
Letting $n \to \infty$, one can get the desired estimates. Note that $a \in L^{\infty}({\bf R})$, so it is non-effective on the norm.
$\Box$
\begin{rem}{\rm Obviously, Proposition \ref{pro2.2} is also true for the weak solution $u(t,x)$ with initial data $[u_{0},u_{1}] \in (H^{1}({\bf R})\cap L^{2}({\bf R},w))\times L^{2}({\bf R},w)$.} 
\end{rem}


\section{An application to semilinear problem }
In this section, we consider the Cauchy problem for semilinear wave equation 
\begin{equation}
	u_{tt}(t,x) - u_{xx}(t,x) + V(x)u(t,x) + a(x)u_{t}(t,x) = |u(t,x)|^p,\ \ \ (t,x)\in (0,\infty)\times {\bf R},\label{eqn2}
\end{equation}
\begin{equation}
	u(0,x)= u_{0}(x),\ \ u_{t}(0,x)= u_{1}(x),\ \ \ x\in {\bf R}.\label{initial2}
\end{equation}
Here, to observe the effect between the potential and the power of nonlinearity, as a trial, we fix the form of $V(x)$ as in Example 1. Indeed, let $V \in C^{1}({\bf R})$ satisfy
\[V(x) = \left\{
  \begin{array}{ll}
   \displaystyle{V_{0}\vert x\vert^{-\beta}},&
       \qquad \vert x\vert \geq L, \\[0.2cm]
   \displaystyle{\frac{2V_{0}}{L^{\beta}}- {\frac{V_{0}}{L^{2\beta}}\vert x\vert^{\beta}}},& \qquad \vert x\vert \leq L,
   \end{array} \right. \]
where $\beta > 1$, $L > 0$ as defined in {\bf (A.2)} and $V_{0} > 0$ is small enough to guarantee the decay estimates of Theorem 1.1.
Under these assumptions on $V(x)$, Theorem 1.1 holds true naturally.

In addition, some important assumptions are imposed on $p > 1$ to derive our main results.\\
({\bf B.1})\, There exists $R>L$ such that 
$$ \text{supp}\,u_0 \cup \text{supp}\, u_1 \subset  B_R:=\{x : |x| \leq R \}.$$
({\bf B.2})\, The exponent $p$ satisfies
\[p>5+2\beta=:p^{*}(\beta).\]
\noindent

\begin{rem}{\rm When $\beta = 2$, the lower bound exponent $p^{*}(2)$ to get the global existence of small data solution is equal to $9$. Note that $\beta = 2$ corresponds to the scale invariant case. As $\beta \to \infty$ (less strong potential), the power $p$ must be chosen large enough.}
\end{rem}
\begin{rem}{\rm In \cite{Ono} a semilinear problem (3.1)-(3.2) with $a(x) =$ constant $ > 0$ and a potential $V(x)$ satisfying $V(x) \geq \frac{k_{0}}{(1+\vert x\vert)^{\lambda}}$ is considered ($k_{0} > 0$). There $\lambda \in [0,\frac{1}{2})$ (long-range potential) and $p \geq 5$ can be treated for $1-D$ case. So, $p^{*}(\beta)$ seems to be reasonable. Note that formally one sees $p^{*}(0) = 5$.}
\end{rem}
\begin{rem}{\rm For each $\beta > 1$, to check a blowup result in the case of $p \in (1,p^{*}(\beta)]$ is still completely open.}
\end{rem}
We prepare useful tools to get the priori estimates of the solution to semilinear problem (3.1)-(3.2).
\begin{lem} [\cite{Ikehata2}, Lemma 2.3] 
If $\theta>1$, there exists a constant $C_\theta>0$ depending only on $\theta$ such that
\begin{equation}\label{28}
	\int_0^t(1+t-s)^{-\frac{1}{2}}(1+s)^{-\theta}ds\leq C_\theta (1+t)^{-\frac{1}{2}}
\end{equation}
for all $t>0$.
\end{lem}

Based on \eqref{28} and the decay estimates obtained for the linear problem \eqref{eqn}-\eqref{initial}, we demonstrate the global existence of small data solution  and decay property for semilinear wave equation. Our main result reads as follows.
\begin{theo}
Let $\beta > 1$. Under the assumptions {\bf (B.1)}, {\bf (B.2)} and Theorem {\rm \ref{th1}}, there exists $\delta>0$ such that if $[u_{0},u_{1}] \in (H^{1}({\bf R})\cap L^{2}({\bf R},w))\times L^{2}({\bf R},w)$ satisfies $I_0< \delta$, the semilinear problem \eqref{eqn2}-\eqref{initial2} admits a unique global solution $u \in {\rm C}([0,+\infty);H^{1}({\bf R}))\cap {\rm C}^{1}([0,+\infty);L^{2}({\bf R}))$ satisfying
\begin{equation}
	\|u(t,\cdot)\|\leq C I_0,
\end{equation}
\begin{equation}
	\|u_t(t,\cdot)\|+\|u_x(t,\cdot)\| +
	\|\sqrt{V(\cdot)}u(t,\cdot)\|
	\leq C I_0(1+t)^{-\frac{1}{2}},
\end{equation}
where 
\[
I_0:=\|u_0\|_{H^1}+\|u_1\|+ \Vert \frac{u_{1}+a(\cdot)u_{0}}{\sqrt{V(\cdot)}}\Vert.\]

\end{theo}
{\it Proof.} 
By a standard semigroup theory, semilinear problem \eqref{eqn2}-\eqref{initial2} can be rewritten as 
\begin{equation} \label{29}
U(t)=S(t)U_0+\int_0^t S(t-s)F(s)ds,
\end{equation}
where $U(t)=[u(t,\cdot), u_t(t,\cdot)]^T$, $U(0)=[u_0, u_1]^T$, $F(s)=[0, |u(s,\cdot)|^p ]^T$ and $S(t)$ denotes the semigroup corresponding to the linear problem. 

For convenience, we introduce the following notation
$$\|U(t)\|_E=\| u_t(t,\cdot)\|+\| u_x(t,\cdot)\|+\|\sqrt{V(\cdot)}u(t,\cdot)\|.$$

It follows from the assumption $p>1$ that there exists a unique mild solution $u \in {\rm C}([0,T);H^{1}({\bf R}))\cap {\rm C}^{1}([0,T);L^{2}({\bf R}))$ for some $T>0$. To show the global existence, it is sufficient to establish the priori estimates for solution and energy in the interval of existence.

We proceed our argument on the bases of \cite{Nakao-1} (see also \cite{Ikehata2}).

Adopting Theorem 1.1, one obtains
\begin{align} \label{27}
	\|U(t)\|_E \leq CI_0(1+t)^{-\frac{1}{2}} +C\int_0^t(1+t-s)^{-\frac{1}{2}}\big(\|\frac{1}{\sqrt{V(\cdot)}}|u(s,\cdot)|^p\| + \| |u(s,\cdot)|^p \| \big).
\end{align}

Next we restrict our attention to the estimates for $\|\frac{1}{\sqrt{V(\cdot)}} |u(s,\cdot)|^p \|$ and $\| |u(s,\cdot)|^p \|$.

By continuity, let us assume that there exists $M > 0$ such that 
\begin{equation} \label{21}
	\|U(t)\|_E \leq M I_{0}(1+t)^{-\frac{1}{2}},\quad t \in [0, T),
\end{equation}
and
\begin{equation}\label{22}
	\|u(t,\cdot)\| \leq M I_{0},\quad t \in [0, T).
\end{equation}
By \eqref{21} and \eqref{22} one does not deny existence of some times $t', t'' \in (0,T)$ such that
\[\|U(t')\|_E = M I_{0}(1+t')^{-\frac{1}{2}}, \quad \|u(t'',\cdot)\| = M I_{0}.\]
Note that \eqref{21} and \eqref{22} are realized if we take large $M>0$ to satisfy 
\[\Vert U(0)\Vert_{E} < MI_{0}, \quad \Vert u_{0}\Vert < MI_{0}.\]
By assumption {$(\bf B.1)$}, one has
$$\text{supp} ~u(s,\cdot) \subset B_{R + s}, ~~~~~~~ s \in [0, T). $$
It follows from the definition of $V(x)$ and the Gagliardo-Nirenberg inequality that
	\begin{align} 
		\|\frac{1}{\sqrt{V(\cdot)}}|u{(s,\cdot)}|^p\|^2 &\leq  V(L)^{-1} \int_{|x| \leq L} |u(s,  x)|^{2p}dx +V_{0}^{-1} \int_{|x| \leq L}  |x|^{\beta}|u(s, x)|^{2p}dx \notag \\
		& \leq 2 \big(V_{0}^{-1}  L^\beta+V_{0}^{-1} (R+s)^{\beta}\big)\int_{R}|u(s, x)|^{2p}dx\notag\\
		& \leq 2 V_{0}^{-1}\big(L^{\beta}+(R+s)^{\beta}\big)\|u(s,\cdot)\|^{2p}_{2p}\notag\\
		& \leq C V_{0}^{-1} \big(L^{\beta}+(R+s)^{\beta}\big) \|u(s,\cdot)\|^{2(1-\theta)p}\|u_{x}(s,\cdot)\|^{2\theta p} \notag, 
	\end{align}
where 
$$\theta=\frac{p-1}{2p}.$$
Therefore
	\begin{align}  \label{23}
	\|\frac{1}{\sqrt{V(\cdot)}}|u{(s,\cdot)}|^p\| \leq C V_{0}^{-\frac{1}{2}} (L^{\frac{\beta}{2}}+(R+s)^{\frac{\beta}{2}}) \|u(s,\cdot)\|^{(1-\theta)p}\|u_{x}(s,\cdot)\|^{\theta p}.
\end{align}
Submitting  \eqref{21} and \eqref{22} to \eqref{23} yields
\begin{align}\label{25}
	\| \frac{1}{\sqrt{V(\cdot)}}|u{(s,\cdot)}|^p \| &\leq C V_{0}^{-\frac{1}{2}}M^{p}I_{0}^{p} \big(L^{\frac{\beta}{2}}+(R+s)^{\frac{\beta}{2}}\big) (1+s)^{-\frac{\theta p}{2}} \notag \\
	& \leq  C V_{0}^{-\frac{1}{2}} M^{p}I_{0}^{p} \big(L^{\frac{\beta}{2}}+(R+s)^{\frac{\beta}{2}}\big) (1+s)^{-\frac{p-1}{4}}.
\end{align}
Similarly, according to \eqref{21}, \eqref{22} and Gagliardo-Nirenberg inequality, one has
\begin{align} \label{26}
\||u(s, \cdot)|^p \|&=
\|u(s, \cdot) \|^p_{2p} \notag \\
& \leq C\|u(s,\cdot)\|^{(1-\theta)p}\|u_{x}(s,\cdot)\|^{\theta p} \notag \\ 
&\leq C M^{p}I_{0}^{p}(1+s)^{-\frac{p-1}{4}}.
\end{align}
Submitting  \eqref{25} and \eqref{26} to \eqref{27} leads to
\begin{align}
	\| U(t) \|_{E} &\leq C I_{0}(1+t)^{-\frac{1}{2}}+C \int_{0}^{t}(1+t-s)^{-\frac{1}{2}}\big(\|u{(s,\cdot)} \|^p_{2p}+ \|\frac{1}{\sqrt{V(\cdot)}} |u{(s,\cdot)}|^{p} \| \big)ds \notag\\
	& \leq C I_{0}(1+t)^{-\frac{1}{2}}+CM^{p}I_{0}^{p}\int_{0}^{t}(1+t-s)^{-\frac{1}{2}}(1+s)^{-\frac{p-1}{4}}ds \notag\\
	& ~~~~+ C V_{0}^{-\frac{1}{2}}M^{p}I_{0}^{p} \int_{0}^{t}(1+t-s)^{-\frac{1}{2}}\big((R+s)^{\frac{\beta}{2}}+L^{\frac{\beta}{2}}\big)(1+s)^{-\frac{p-1}{4}} ds \notag\\
	& \leq C I_{0}(1+t)^{-\frac{1}{2}}+CM^{p}I_{0}^{p}(V_0^{-\frac{1}{2}}R^{\frac{\beta}{2}}+1)\int_{0}^{t}(1+t-s)^{-\frac{1}{2}}(1+s)^{-\frac{p-1}{4}+\frac{\beta}{2}}ds
\end{align}
for all $t \in [0,T)$.

Let 
$$\gamma :=\frac{p-1}{4}-\frac{\beta}{2}.$$
By assumption {$(\bf B.2)$}, we see $\gamma>1$. Using \eqref{28} one has
\begin{align}
	\|U(t) \|_{E} \leq &\big(C I_{0}+CM^{p}I_{0}^{p}(V_0^{-\frac{1}{2}}R^{\frac{\beta}{2}}+1)\big)(1+t)^{-\frac{1}{2}}  \notag\\
	 = &C I_{0}\big(1+ M^{p}I_{0}^{p-1}(V_0^{-\frac{1}{2}}R^{\frac{\beta}{2}}+1) \big)(1+t)^{-\frac{1}{2}} \notag\\
	 =& I_{0}Q_{0}(I_{0},M,R,V_0)(1+t)^{-\frac{1}{2}}
\end{align}
for all $t \in [0,T)$, where 
$$Q_{0}(I_{0},M,R,V_0)=C \big(1+ M^{p}I_{0}^{p-1}(V_0^{-\frac{1}{2}}R^{\frac{\beta}{2}}+1) \big).$$
\\
Next, we derive $L^{2}$-bound for the local solution to problem \eqref{eqn2}-\eqref{initial2}.

In fact, one has from Theorem 1.1 and \eqref{29} that
\begin{equation} \label{30}
		\| u(t,\cdot) \| \leq C I_{0}+C \int_{0}^{t}\big(\|u{(s,\cdot)} \|^p_{2p}+ \|\frac{1}{\sqrt{V(\cdot)}} |u{(s,\cdot)}|^{p} \| \big)ds.
\end{equation}
Substituting \eqref{25} and \eqref{26} to \eqref{30} and proceeding similar arguments to the estimates for $\|U(t)\|_E$, one has
\begin{align*}
	\|u(t,\cdot)\|
	\leq I_{0}Q_{0}(I_{0},M,R,V_0)
\end{align*}
for all $t \in [0,T)$.

Take $ M > 0$ further to satisfy $M > C$, and $I_{0} $ small enough such that
\begin{align}\label{33}
	CM^{p}I_0^{p-1}(V_0^{-\frac{1}{2}}R^{\frac{\beta}{2}}+1) < M - C.
\end{align}
By choosing 
$$\delta := \left(\frac{M-C}{CM^{p}(V_0^{-\frac{1}{2}}R^{\frac{\beta}{2}}+1)}\right)^{\frac{1}{p-1}},$$
and if $I_0 < \delta$, then it follows from \eqref{33} that 
$$Q_{0}(I_{0},M,R,V_0) < M.$$
Therefore, we see that
\begin{equation}\label{31}
\Vert U(t)\Vert_{E} < M I_{0}(1+t)^{-\frac{1}{2}}\quad \text{in}\quad [0,T), \\
\end{equation}	
\begin{equation}\label{32}
\|u(t,\cdot)\| < M I_{0}\quad\text{in}\quad[0,T). 
\end{equation}
This contradicts \eqref{21}-\eqref{22}, and so \eqref{31} and \eqref{32} are true in $[0,T)$. This shows that the local solution can be extended globally in time and the estimates \eqref{31} and \eqref{32} for the solution $u(t,x)$ hold true for all $t \geq 0$, which completes the proof.
\hfill
$\Box$


\vspace{0.5cm}
\noindent{\em Acknowledgement.}
\smallskip
The work of the first author (R. IKEHATA) was supported in part by Grant-in-Aid for Scientific Research (C) 22540193 of JSPS.


\end{document}